# ONE-STEP SPARSE ESTIMATES IN NONCONCAVE PENALIZED LIKELIHOOD MODELS[1]

By Hui Zou and Runze Li

*University of Minnesota and Pennsylvania State University*

Fan and Li propose a family of variable selection methods via penalized likelihood using concave penalty functions. The nonconcave penalized likelihood estimators enjoy the oracle properties, but maximizing the penalized likelihood function is computationally challenging, because the objective function is nondifferentiable and nonconcave. In this article, we propose a new unified algorithm based on the local linear approximation (LLA) for maximizing the penalized likelihood for a broad class of concave penalty functions. Convergence and other theoretical properties of the LLA algorithm are established. A distinguished feature of the LLA algorithm is that at each LLA step, the LLA estimator can naturally adopt a sparse representation. Thus, we suggest using the one-step LLA estimator from the LLA algorithm as the final estimates. Statistically, we show that if the regularization parameter is appropriately chosen, the one-step LLA estimates enjoy the oracle properties with good initial estimators. Computationally, the one-step LLA estimation methods dramatically reduce the computational cost in maximizing the nonconcave penalized likelihood. We conduct some Monte Carlo simulation to assess the finite sample performance of the one-step sparse estimation methods. The results are very encouraging.

**1. Introduction.** Variable selection and feature extraction are fundamental for knowledge discovery and predictive modeling with high-dimensionality (Fan and Li [13]). The best subset selection procedure along with traditional model selection criteria, such as AIC and BIC, becomes infeasible for feature

Received May 2006; revised March 2007.

[1]Supported in part by NSF Grant DMS-03-48869 and partially supported by a National Institute on Drug Abuse (NIDA) Grant P50 DA10075.

[2]Discussed in 10.1214/07-AOS0316A, 10.1214/07-AOS0316B and 10.1214/07-AOS0316C; rejoinder at 10.1214/07-AOS0316REJ.

selection from high-dimensional data due to too expensive computational cost. Furthermore, the best subset selection suffers from several drawbacks, the most severe of which is its lack of stability as analyzed in Breiman [4]. LASSO (Tibshirani [32]) method utilizes the $L_1$ penalty to automatically select significant variable via continuous shrinkage, thus retaining the good features of both the best subset selection and ridge regression. In the same spirit of LASSO, the penalized likelihood with nonconcave penalty functions has been proposed to select significant variables for various parametric models, including generalized linear regression models and robust linear regression model (Fan and Li [10] and Fan and Peng [15]), and some semiparametric models, such as the Cox model and partially linear models (Fan and Li [11, 12] and Cai, Fan, Li and Zhou [5]). Fan and Li [10] provide deep insights into how to select a penalty function. They further advocate the use of penalty functions satisfying certain mathematical conditions such that the resulting penalized likelihood estimate possesses the properties of sparsity, continuity and unbiasedness. These mathematical conditions imply that the penalty function has to be singular at the origin and nonconvex over $(0, \infty)$. In the work aforementioned, it has been shown that when the regularization parameter is appropriately chosen, the nonconcave penalized likelihood estimates perform as well as the oracle procedure in terms of selecting the correct subset model and estimating the true nonzero coefficients.

Although nonconcave penalized likelihood approaches have promising theoretical properties, the singularity and nonconvexity of the penalty function challenge us to invent numerical algorithms which are capable of maximizing a nondifferentiable nonconcave function. Fan and Li [10] suggested iteratively, locally approximating the penalty function by a quadratic function and referred such approximation as to local quadratic approximation (LQA). With the aid of the LQA, the optimization of penalized likelihood function can be carried out using a modified Newton–Raphson algorithm. However, as pointed out in Fan and Li [10] and Hunter and Li [20], the LQA algorithm shares a drawback of backward stepwise variable selection: If a covariate is deleted at any step in the LQA algorithm, it will necessarily be excluded from the final selected model (see Section 2.2 for more details). Hunter and Li [20] addressed this issue by optimizing a slightly perturbed version of LQA, which alleviates the aforementioned drawback, but it is difficult to choose the size of perturbation. Another strategy to overcome the computational difficulty is using the one-step (or $k$-step) estimates from the iterative LQA algorithm with good starting estimators, as suggested by Fan and Li [10]. This is similar to the well-known one-step estimation argument in the maximum likelihood estimation (MLE) setting (Bickel [2], Lehmann and Casella [24], Robinson [30] and Cai, Fan, Zhou and Zhou [6]). See also Fan and Chen [9], Fan, Lin and Zhou [14] and Cai et al. [6] for some recent work on one-step estimators in local and marginal likelihood models. However, the problem with the one-step LQA estimator is that it cannot have a sparse representation, thus losing the most attractive and important property of the nonconcave penalized likelihood estimator.

In this article we develop a methodology and theory for constructing an efficient one-step sparse estimation procedure in nonconcave penalized likelihood models. For that purpose, we first propose a new iterative algorithm based on local linear approximation (LLA) for maximizing the nonconcave penalized likelihood. The LLA enjoys three significant advantages over the LQA and the perturbed LQA. First, in the LLA we do not have to delete any small coefficient or choose the size of perturbation in order to avoid numerical instability. Second, we demonstrate that the LLA is *the best convex* minorization–maximization (MM) algorithm, thus proving the convergence of the LLA algorithm by the *ascent property* of MM algorithms (Lange, Hunter and Yang [23]). Third, the LLA naturally produces a sparse estimates via continuous penalization. We then propose using the one-step LLA estimator from the LLA algorithm as the final estimates. Computationally, the one-step LLA estimates alleviate the computation burden in the iterative algorithm and overcome the potential local maxima problem in maximizing the nonconcave penalized likelihood. In addition, we can take advantage of the efficient algorithm for solving LASSO to compute the one-step LLA estimator. Statistically, we show that if the regularization parameter is appropriately chosen, the one-step LLA estimates enjoy the oracle properties, provided that the initial estimates are good enough. Therefore, the one-step LLA estimator can dramatically reduce the computation cost without losing statistical efficiency.

The rest of the paper is organized as follows. In Section 2 we introduce the local linear approximation algorithm and discuss its various properties. In Section 3 we discuss the one-step LLA estimator, in which asymptotical normality and consistency of selection are established. Section 4 describes the implementation detail, and Section 5 shows numerical examples. Proofs are presented in Section 6.

**2. Local linear approximation algorithm.** Suppose that $\{(\mathbf{x}_i, y_i)_{i=1}^n\}$ are $n$ identically and independently distributed samples, where $\mathbf{x}_i$ denotes the $p$-dimension predictor and $y_i$ is the response variable. Assume that $y_i$ depends on $\mathbf{x}_i$ through a linear combination $\mathbf{x}_i^\mathsf{T}\boldsymbol{\beta}$, and the conditional log-likelihood given $\mathbf{x}_i$ is $\ell_i(\boldsymbol{\beta}, \boldsymbol{\phi}) = \ell_i(\mathbf{x}_i^\mathsf{T}\boldsymbol{\beta}, y_i, \boldsymbol{\phi})$, where $\boldsymbol{\phi}$ is a dispersion parameter. In some models, such as logistic regression and Poisson regression, there is no dispersion parameter. In linear regression model, $\boldsymbol{\phi}$ is the variance of the random error, and is often estimated separately after $\boldsymbol{\beta}$ is estimated. In most variable selection applications, we do not penalize the dispersion parameter (Frank and Friedman [16], Tibshirani [32], Fan and Li [10] and Miller [28]). Thus, we simplify notation in the reminder of this paper by suppressing $\boldsymbol{\phi}$, and further use $\ell_i(\boldsymbol{\beta})$ to stand for $\ell_i(\mathbf{x}_i^\mathsf{T}\boldsymbol{\beta}, y_i, \boldsymbol{\phi})$.

2.1. *Penalized likelihood.* In the variable selection problem, the assumption is that some components of $\boldsymbol{\beta}$ are zero. The goal is to identify and estimate the subset model. In this work, we consider the variable selection methods by maximizing the penalized likelihood function taking the form

$$Q(\boldsymbol{\beta}) = \sum_{i=1}^{n} \ell_i(\boldsymbol{\beta}) - n \sum_{j=1}^{p} p_{\lambda_j}(|\beta_j|). \tag{2.1}$$

In principle, $p_{\lambda_j}$ can be different for different components (coefficients). For ease of presentation, we let $p_{\lambda_j}(|\beta_j|) = p_\lambda(|\beta_j|)$, that is, the same penalty function is applied to every component of $\boldsymbol{\beta}$. Formulation in (2.1) includes

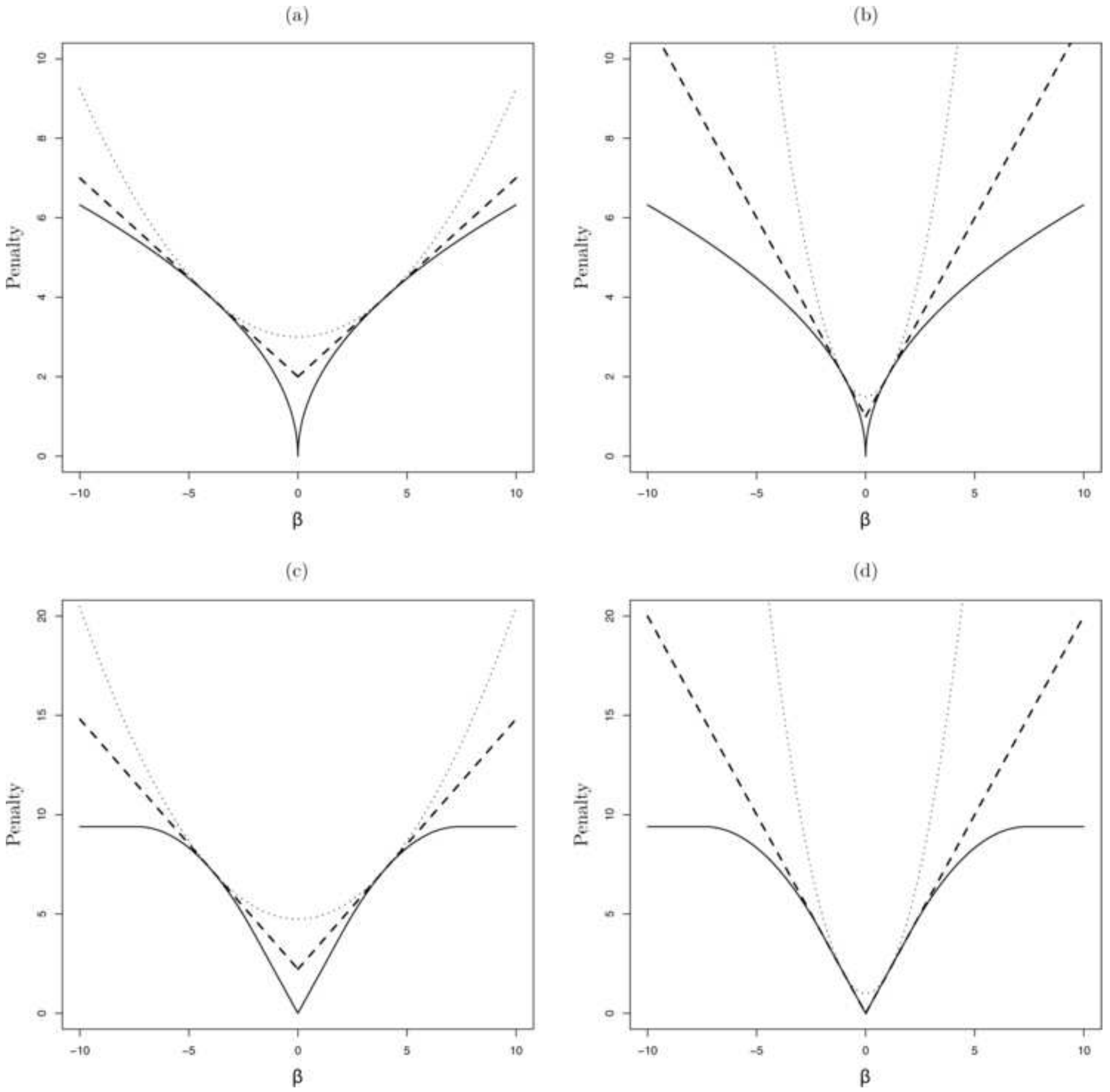

FIG. 1. *Plot of local quadratic approximation (thin dotted lines) and local linear approximation (thick broken lines) at $\beta = 4$ and 1.* (a) *and* (b) *are for the $L_{0.5}$ penalty with $\lambda = 2$, and* (c) *and* (d) *are for the SCAD penalty with $\lambda = 2$.*

many popular variable selection methods. For instance, the best subset selection amounts to using the $L_0$ penalty, while the LASSO (Tibshirani [32]) uses the $L_1$ penalty $p_\lambda(|\beta|) = \lambda|\beta|$. Bridge regression (Frank and Friedman [16]) uses the $L_q$ penalty $p_\lambda(|\beta|) = \lambda|\beta|^q$. When $0 < q < 1$, the $L_q$ penalty is concave over $(0, \infty)$, and nondifferentiable at zero. The SCAD penalty (Fan and Li [10]) is a concave function defined by $p_\lambda(0) = 0$ and for $|\beta| > 0$

$$p'_\lambda(|\beta|) = \lambda I(|\beta| \le \lambda) + \frac{(a\lambda - |\beta|)_+}{a-1} I(|\beta| > \lambda) \qquad \text{for some } a > 2. \tag{2.2}$$

Often $a = 3.7$ is used. The notation $z_+$ stands for the positive part of $z$: $z_+$ is $z$ if $z > 0$, zero otherwise. The SCAD penalty and $L_{0.5}$ penalty are illustrated in Figure 1. Note that with a concave penalty the penalized likelihood in (2.1) is a nonconcave function. Hence maximizing nonconcave penalized likelihood is challenging. Antoniadis and Fan [1] proposed nonlinear regularized Sobolev interpolators (NRSI) and regularized one-step estimator (ROSE) for nonconvex penalized least squares problems under wavelets settings. They further introduced the graduated nonconvexity (GNC) algorithm for minimizing high-dimensional nonconvex penalized least squares problem. The GNC algorithm was first developed for reconstructing piecewise continuous images (Black and Zisserman [3]). The GNC algorithm offers nice ideas for minimizing high-dimensional nonconvex objective function, but in general, it is computationally intensive, and its implementation depends on a sequences of tuning parameters. Fan and Li [10] proposed the local quadratic approximation (LQA) algorithm for the nonconcave penalized likelihood. We introduce the LQA algorithm in Section 2.2 in detail. Hunter and Li [20] showed that the LQA shares the same spirit as that of the MM algorithm (Lange et al. [23]). Wu [33] pointed out that the MM algorithm and GNC algorithm share the same spirit in terms of optimization transfer. In general, the GNC algorithms do not guarantee the ascent property for maximization problems, evidenced from Figure 8(c) in Antoniadis and Fan [1], while the MM algorithms enjoy the ascent property, as demonstrated in Hunter and Li [20].

2.2. *Local quadratic approximation.* It can be seen from Figure 1 that the penalized likelihood functions become nondifferentiable at the origin and nonconcave with respect to $\boldsymbol{\beta}$. The singularity and nonconcavity make it difficult to maximize the penalized likelihood functions. Suppose that we are given an initial value $\boldsymbol{\beta}^{(0)}$ that is close to the true value of $\boldsymbol{\beta}$. Fan and Li (2001) propose locally approximating the first order derivative of the penalty function by a linear function:

$$[p_\lambda(|\beta_j|)]' = p'_\lambda(|\beta_j|)\,\mathrm{sign}(\beta_j) \approx \{p'_\lambda(|\beta_j^{(0)}|)/|\beta_j^{(0)}|\}\beta_j.$$

Thus, they use a LQA to the penalty function:

$$p_\lambda(|\beta_j|) \approx p_\lambda(|\beta_j^{(0)}|) + \tfrac{1}{2}\{p'_\lambda(|\beta_j^{(0)}|)/|\beta_j^{(0)}|\}(\beta_j^2 - \beta_j^{(0)2}) \quad \text{for } \beta_j \approx \beta_j^{(0)}. \tag{2.3}$$

Figure 1 illustrates the LQA for the $L_{0.5}$ penalty and the SCAD penalty. With iteratively updating the LQA, Newton–Raphson algorithm can be modified for maximization of the penalized likelihood function. Specifically, we take the unpenalized likelihood estimate to be the initial value $\boldsymbol{\beta}^{(0)}$: For $k = 1, 2, \ldots,$ repeatedly solve

$$\boldsymbol{\beta}^{(k+1)} = \arg\max\left\{\sum_{i=1}^n \ell_i(\boldsymbol{\beta}) - n\sum_{j=1}^p \frac{p'_\lambda(|\beta_j^{(k)}|)}{2|\beta_j^{(k)}|}\beta_j^2\right\}. \tag{2.4}$$

Stop the iteration if the sequence of $\{\boldsymbol{\beta}^{(k)}\}$ converges.

To avoid numerical instability, Fan and Li [10] suggested that if $\beta_j^{(k)}$ in (2.4) is very close to 0, say $|\beta_j^{(k)}| < \varepsilon_0$ (a prespecified value), then set $\hat{\beta}_j = 0$ and delete the $j$th component of $\mathbf{x}$ from the iteration. Thus, the LQA algorithm shares a drawback of backward stepwise variable selection: if a covariate is deleted at any step in the LQA algorithm, it will necessarily be excluded from the final selected model. Furthermore, one has to choose $\varepsilon_0$, which practically becomes an additional tuning parameter. The size of $\varepsilon_0$ potentially affects the degree of sparsity of the solution as well as the speed of convergence. Hunter and Li [20] studied the convergence property of the LQA algorithm. They found that the LQA algorithm is one of minorize–maximize (MM) algorithms, extensions of the well-known EM algorithm. They further demonstrated that the behavior of the LQA algorithm is the same as that of an EM algorithm with the LQA playing the same role of E-step in the EM algorithm. To avoid numerical instability and the drawback of backward stepwise variable selection, Hunter and Li [20] suggested optimizing a slightly perturbed version of (2.4) bounding the denominator away from zero: for $k = 1, 2, \ldots,$ repeatedly solve

$$\boldsymbol{\beta}^{(k+1)} = \arg\max\left\{\sum_{i=1}^n \ell_i(\boldsymbol{\beta}) - n\sum_{j=1}^p \frac{p'_\lambda(|\beta_j^{(k)}|)}{2\{|\beta_j^{(k)}| + \tau_0\}}\beta_j^2\right\}, \tag{2.5}$$

for a prespecified size perturbation $\tau_0$. Stop the iteration if the sequence of $\{\boldsymbol{\beta}^{(k)}\}$ converges. In the practical implementation, we have to determine the size of perturbation. This sometimes may be difficult, and furthermore, the size of $\tau_0$ potentially affects the degree of sparsity of the solution as well as the speed of convergence.

2.3. *Local linear approximation.* To eliminate the weakness of the LQA, we propose a new unified algorithm based on local linear approximation to the penalty function:

$$p_\lambda(|\beta_j|) \approx p_\lambda(|\beta_j^{(0)}|) + p'_\lambda(|\beta_j^{(0)}|)(|\beta_j| - |\beta_j^{(0)}|), \qquad \text{for } \beta_j \approx \beta_j^{(0)}. \tag{2.6}$$

Figure 1 illustrates the LLA for the $L_{0.5}$ penalty and the SCAD penalty. Fan and Li [10] show that in order to have a continuous thresholding rule, the penalty function must satisfy a *continuity condition*: the minimum of $|\theta| + p'_\lambda(|\theta|)$ is attained at zero. Although the $L_{0.5}$ penalty fails to hold the continuity condition, we show in Section 3 that it is still good for deriving continuous one-step sparse estimates. For ease of presentation, we assume in this section, unless otherwise specified, that the right derivative of $p_\lambda(\cdot)$ at 0 is finite.

Similar to the LQA algorithm, the maximization of the penalized likelihood can be carried out as follows. Set the initial value $\boldsymbol{\beta}^{(0)}$ be the unpenalized maximum likelihood estimate. For $k = 1, 2, \ldots$, repeatedly solve

$$\boldsymbol{\beta}^{(k+1)} = \arg\max \left\{ \sum_{i=1}^{n} \ell_i(\boldsymbol{\beta}) - n \sum_{j=1}^{p} p'_\lambda(|\beta_j^{(k)}|)|\beta_j| \right\}. \tag{2.7}$$

Stop the iterations if the sequence of $\{\boldsymbol{\beta}^{(k)}\}$ converges. We refer this algorithm as to the LLA algorithm. The LLA algorithm is distinguished from the LQA algorithm in that $\boldsymbol{\beta}^{(k+1)}$ and the final estimates naturally adopt a sparse representation. The LLA algorithm inherits the good features of LASSO in terms of computational efficiency, and therefore the maximization can be solved by efficient algorithms, such as the least angle regression (LARS) algorithm (Efron, Hastie, Johnstone and Tibshirani [8]). From (2.7), the approximation is numerical stable, and thus, the drawback of backward variable selection can be avoided in LLA algorithm.

We next study the convergence of the LLA algorithm. Denote

$$\phi^*(\beta_j|\beta_j^{(0)}) = p_\lambda(|\beta_j^{(0)}|) + p'_\lambda(|\beta_j^{(0)}|)(|\beta_j| - |\beta_j^{(0)}|)$$

and

$$G(\boldsymbol{\beta}|\boldsymbol{\beta}^{(k)}) = \sum_{i=1}^{n} \ell_i(\boldsymbol{\beta}) - n \sum_{j=1}^{p} \phi^*(\beta_j|\beta_j^{(k)}).$$

THEOREM 1. *For a differentiable concave penalty function $p_\lambda(\cdot)$ on $[0, \infty)$, we have*

$$Q(\boldsymbol{\beta}) \geq G(\boldsymbol{\beta}|\boldsymbol{\beta}^{(k)}) \quad \textit{and} \quad Q(\boldsymbol{\beta}^{(k)}) = G(\boldsymbol{\beta}^{(k)}|\boldsymbol{\beta}^{(k)}). \tag{2.8}$$

*Furthermore, the LLA has the ascent property, that is, for all* $k = 0, 1, 2, \ldots$

$$Q(\boldsymbol{\beta}^{(k+1)}) \geq Q(\boldsymbol{\beta}^{(k)}). \tag{2.9}$$

*If the penalty function is strictly concave then we always take* ">" *in (2.9).*

From (2.8), $G(\boldsymbol{\beta}|\boldsymbol{\beta}^{(k)})$ is a minorization of $Q(\boldsymbol{\beta})$, and finding $\boldsymbol{\beta}^{(k+1)}$ is the maximize-step in MM (minorize–maximize) algorithms. Therefore, the LLA algorithm is an instance of the MM algorithms. For a survey of work in MM algorithms, see Heiser [19] and Lange, Hunter and Yang [23].

The analysis of convergence of LLA can be done by following the general convergence results for MM algorithms. Let $M(\boldsymbol{\beta})$ denote the map defined by the LLA algorithm from $\boldsymbol{\beta}^{(k)}$ to $\boldsymbol{\beta}^{(k+1)}$. Note that the penalty function has continuous first derivative and solving $\boldsymbol{\beta}^{(k+1)}$ is a convex optimization problem, thus $M$ is a continuous map. We define a stationary point of the function $Q(\boldsymbol{\beta})$ to be any point $\boldsymbol{\beta}$ at which the gradient vectors is zero.

PROPOSITION 1. *Given an initial value* $\boldsymbol{\beta}^{(0)}$, *let* $\boldsymbol{\beta}^{(k)} = M^k(\boldsymbol{\beta}^{(0)})$. *If* $Q(\boldsymbol{\beta}) = Q(M(\boldsymbol{\beta}))$ *only for stationary points of* $Q$ *and if* $\boldsymbol{\beta}^*$ *is a limit point of the sequence* $\{\boldsymbol{\beta}^{(k)}\}$, *then* $\boldsymbol{\beta}^*$ *is a stationary point of* $Q(\boldsymbol{\beta})$.

Proposition 1 is a slightly modified version of Lyapunov's theorem in Lange [22]. We omit its proof. In Theorem 1, we show that the LLA of $p_\lambda(\cdot)$ provides a majorization of the penalty function $p_\lambda(\cdot)$. In fact, the LLA is the best convex majorization of $p_\lambda(\cdot)$ as stated in the next theorem.

THEOREM 2. *Denote by* $\psi^*(\cdot)$ *the LLA approximation of* $p_\lambda(\cdot)$. $\psi^*(t) = p_\lambda(t_0) + p'_\lambda(t_0)(t - t_0)$, $t, t_0 \geq 0$. *Suppose that* $\psi(\cdot)$ *is a convex majorization function of* $p_\lambda(\cdot)$ *at* $t_0$, *that is,*

$$\psi(t_0) = p_\lambda(t_0) \quad \textit{and} \quad \psi(t) \geq p_\lambda(t) \qquad \textit{for all } t.$$

*We must have* $\psi(t) \geq \psi^*(t)$ *for all* $t$. *If the right derivative of* $p_\lambda(\cdot)$ *at zero diverges, the above conclusions hold for* $t_0 > 0$ *and* $t \geq 0$.

Figure 1 shows an illustration of Theorem 2 with the SCAD and $L_{0.5}$ penalties. As can be seen from Figure 1, the LLA approximation is underneath the LQA approximation in all four cases.

The ascent property of MM indicates that MM is an extension of the famous EM algorithm. Under certain conditions, we show that the LLA algorithm can be cast as an EM algorithm.

Suppose that $\exp(-np_\lambda(\cdot))$ is a Laplace transformation of some nonnegative function $H(\cdot)$. Then $H(\cdot)$ is the inverse Laplace transformation of $\exp(-np_\lambda(\cdot))$ and

$$\exp(-np_\lambda(|\beta|)) = \int_0^\infty H(t) e^{-t|\beta|}\, dt. \tag{2.10}$$

For example, if $p_\lambda(|\beta|) = \lambda|\beta|^q$, the Bridge penalty $(0 < q < 1)$, then

$$\exp(-n\lambda|\beta|^q) = \int_0^\infty H(t)e^{-t|\beta|}\,dt,$$

where $H(t) \propto \frac{(n\lambda)^{1/q}}{2} S(\frac{(n\lambda)^{1/q}}{2}t)$ and $S(\cdot)$ is the density of the stable distribution of index $q$ (Mike [27]).

Let $\pi(t) = \frac{2}{t}H(\frac{1}{t})$ and we independently put a Laplacian prior on $\beta_j$

$$p(\beta_j|\tau_j) = \frac{1}{2\tau_j}e^{-|\beta_j|/\tau_j}. \tag{2.11}$$

Further regard $\pi$ as a hyper-prior on $\tau_j$. Then (2.10) implies

$$\exp(-np_\lambda(|\beta_j|)) = \int_0^\infty p(\beta_j|\tau_j)\pi(\tau_j)\,d\tau_j. \tag{2.12}$$

Maximizing $Q(\boldsymbol{\beta})$ is equivalent to computing the posterior mode of $p(\boldsymbol{\beta}|\mathbf{y})$, if we treat $\exp(-np_\lambda(|\beta_j|))$ as the marginal prior of $\beta$. The identity (2.12) implies an EM algorithm for maximizing the posterior $p(\boldsymbol{\beta}|\mathbf{y})$.

To derive the EM algorithm, we consider $\tau_1, \ldots, \tau_p$ as missing data. The complete log-likelihood function (CLF) is

$$\sum_{i=1}^n \ell_i(\boldsymbol{\beta}) + \sum_{j=1}^p \biggl[-\log(2\tau_j) - \frac{|\beta_j|}{2\tau_j} + \log\pi(\tau_j)\biggr].$$

Suppose the current estimator is $\boldsymbol{\beta}^{(k)}$. The E-step computes the conditional mean of CLF

$$E_{\tau|\beta^{(k)},\mathbf{y}}[\mathrm{CLF}] = \sum_{i=1}^n \ell_i(\boldsymbol{\beta}) + \sum_{j=1}^p E\biggl[-\log(2\tau_j) - \frac{|\beta_j|}{\tau_j} + \log\pi(\tau_j)|\beta^{(k)},\mathbf{y}\biggr].$$

The M-step finds $\boldsymbol{\beta}^{(k+1)}$ maximizing $E_{\tau|\beta^{(k)},\mathbf{y}}[\mathrm{CLF}]$. Thus

$$\boldsymbol{\beta}^{(k+1)} = \arg\max \sum_{i=1}^n \ell_i(\boldsymbol{\beta}) + \sum_{j=1}^p \biggl(-|\beta_j| E\biggl[\frac{1}{\tau_j}\bigg|\beta^{(k)},\mathbf{y}\biggr]\biggr). \tag{2.13}$$

THEOREM 3. *Suppose that (2.10)–(2.13) hold for $p_\lambda(\cdot)$, the LLA algorithm and the EM algorithm are identical. Moreover, (2.10) implies that $p_\lambda(\cdot)$ must be a strictly increasing function on $[0, \infty)$ and unbounded. Thus the SCAD penalty does not have an inverse Laplace transformation.*

In the above discussion, we have assumed all the necessary conditions to ensure the the EM algorithm is proper. If this is the case, then Theorem 3 shows that the EM algorithm is exactly the LLA algorithm. On the other hand, it is also worth noting that there are concave penalty functions for which (2.10) cannot be true. The SCAD penalty is such an example. Thus, Theorem 3 also indicates that MM algorithms are more flexible than EM algorithms.

**3. One-step sparse estimates.** In this section, we propose the one-step LLA estimator, which is significantly distinguished from the one-step or $k$-step LQA estimate because it automatically adopts a sparse representation. Thus it can be used as a model selector. One may further define $k$-step LLA estimator, but, in general, it is unnecessary. As demonstrated in Fan and Chen [9] and Cai, Fan and Li [7], both empirically and theoretically, the one-step method is as efficient as the fully iterative method, provided that the initial estimators are reasonably good. In LQA finding $\boldsymbol{\beta}^{(k+1)}$ is a ridge regression problem, which indicates that almost surely, none of the components of $\boldsymbol{\beta}^{(k+1)}$ will be *exact zero*. Hence the one-step or $k$-step LQA estimates in the LQA will not be able to achieve the goal of variable selection. To get insights into the one-step LLA estimator, let us start with linear regression models and consider the penalized least squares.

3.1. *Linear regression models.* The LLA algorithm naturally provides a sparse one-step estimator. For simplicity, let the initial estimate $\boldsymbol{\beta}^{(0)}$ be ordinary least squares estimator. Then the one-step estimator is obtained by

$$\boldsymbol{\beta}^{(1)} = \arg\min \tfrac{1}{2}\|\mathbf{y} - \mathbf{X}\boldsymbol{\beta}\|^2 + n\sum_{j=1}^{p} p'_\lambda(|\beta_j^{(0)}|)|\beta_j|. \tag{3.1}$$

We denote by $\widehat{\boldsymbol{\beta}}(\mathrm{ose})$ the one-step estimator $\boldsymbol{\beta}^{(1)}$.

We show that the one-step estimator enjoys the oracle properties. To this end we assume two regularity conditions:

(A1). $y_i = \mathbf{x}_i\boldsymbol{\beta}_0 + \epsilon_i$, where $\epsilon_1, \ldots, \epsilon_n$ are independent and identically distributed random variables with mean 0 and variance $\sigma^2$,

(A2). $\frac{1}{n}\mathbf{X}^\mathsf{T}\mathbf{X} \to \mathbf{C}$ where $\mathbf{C}$ is a positive definite matrix.

Without loss of generality, let $\boldsymbol{\beta}_0 = (\beta_{01}, \ldots, \beta_{0p})^\mathsf{T} = (\boldsymbol{\beta}_{10}^\mathsf{T}, \boldsymbol{\beta}_{20}^\mathsf{T})^\mathsf{T}$ and $\boldsymbol{\beta}_{20} = 0$. We write

$$\mathbf{C} = \begin{bmatrix} \mathbf{C}_{11} & \mathbf{C}_{12} \\ \mathbf{C}_{21} & \mathbf{C}_{22} \end{bmatrix}.$$

THEOREM 4. *Let $p_{\lambda_n}(\cdot)$ be the SCAD penalty. If $\sqrt{n}\lambda_n \to \infty$ and $\lambda_n \to 0$, then the one-step SCAD estimates $\widehat{\boldsymbol{\beta}}(\mathrm{ose})$ must satisfy:*

(a) *Sparsity: with probability tending to one,* $\widehat{\boldsymbol{\beta}}(\text{ose})_2 = 0$.
(b) *Asymptotic normality:* $\sqrt{n}(\widehat{\boldsymbol{\beta}}(\text{ose})_1 - \boldsymbol{\beta}_{10}) \to N(0, \sigma^2 \mathbf{C}_{11}^{-1})$.

*In addition, consider* $p_{\lambda_n}(\cdot) = \lambda_n p(\cdot)$. *Suppose* $p'(\cdot)$ *is continuous on* $(0, \infty)$ *and there is some* $s > 0$ *such that* $p'(\theta) = O(\theta^{-s})$ *as* $\theta \to 0+$. *Then* (a) *and* (b) *hold, if* $n^{(1+s)/2}\lambda_n \to \infty$ *and* $\sqrt{n}\lambda_n \to 0$.

3.2. *Penalized likelihood.* For a general likelihood model, let $\ell(\boldsymbol{\beta}) = \sum_{i=1}^n \ell_i(\boldsymbol{\beta})$ denote the log-likelihood. Suppose that the log-likelihood function is smooth and has the first two derivatives with respect to $\boldsymbol{\beta}$. For a given initial value $\boldsymbol{\beta}^{(0)}$, the log-likelihood function can be locally approximated by

$$\begin{aligned}\ell(\boldsymbol{\beta}) &\approx \ell(\boldsymbol{\beta}^{(0)}) + \nabla\ell(\boldsymbol{\beta}^{(0)})^{\mathsf{T}}(\boldsymbol{\beta} - \boldsymbol{\beta}^{(0)}) \\ &\quad + \tfrac{1}{2}(\boldsymbol{\beta} - \boldsymbol{\beta}^{(0)})^{\mathsf{T}}\nabla^2\ell(\boldsymbol{\beta}^{(0)})(\boldsymbol{\beta} - \boldsymbol{\beta}^{(0)}).\end{aligned} \tag{3.2}$$

Let us take $\boldsymbol{\beta}^{(0)} = \widehat{\boldsymbol{\beta}}(\text{mle})$. Then $\nabla\ell(\boldsymbol{\beta}^{(0)}) = 0$ by the definition of MLE. Thus, $\boldsymbol{\beta}^{(1)}$ is given by

$$\begin{aligned}\boldsymbol{\beta}^{(1)} = \arg\min \tfrac{1}{2}&(\boldsymbol{\beta} - \boldsymbol{\beta}^{(0)})^{\mathsf{T}}[-\nabla^2\ell(\boldsymbol{\beta}^{(0)})](\boldsymbol{\beta} - \boldsymbol{\beta}^{(0)}) \\ &+ n\sum_{j=1}^p p'_\lambda(|\beta_j^{(0)}|)|\beta_j|.\end{aligned} \tag{3.3}$$

It is interesting to see that (3.3) reduces to the one-step estimates in linear regression models, if we are willing to assume that $\epsilon \sim N(0, \sigma^2)$. However, it should be noted that normality assumption is not needed in Theorem 4.

We show that in the general likelihood setting, $\boldsymbol{\beta}^{(1)}$ is desired the one-step estimates, denoted by $\widehat{\boldsymbol{\beta}}(\text{ose})$. Let $I(\boldsymbol{\beta}_0)$ be the Fisher information matrix and $I_1(\boldsymbol{\beta}_{10}) = I_1(\boldsymbol{\beta}_{10}, 0)$ denote the Fisher information knowing $\boldsymbol{\beta}_{20} = 0$. Note that $I(\boldsymbol{\beta}_0)$ is a $p \times p$ matrix and $I_1(\boldsymbol{\beta}_{10})$ is a submatrix of $I(\boldsymbol{\beta}_0)$. It is well known that under some regularity conditions (Lehmann and Casella [24]), $n^{-1}\nabla^2\ell(\widehat{\boldsymbol{\beta}}(\text{mle})) \to_P -I(\boldsymbol{\beta}_0)$, and

$$\sqrt{n}(\boldsymbol{\beta}_0 - \widehat{\boldsymbol{\beta}}(\text{mle})) \xrightarrow{D} W = N(0, I^{-1}(\boldsymbol{\beta}_0)).$$

THEOREM 5. *Let* $p_{\lambda_n}(\cdot)$ *be the SCAD penalty. If* $\sqrt{n}\lambda_n \to \infty$ *and* $\lambda_n \to 0$, *then the one-step SCAD estimates* $\widehat{\boldsymbol{\beta}}(\text{ose})$ *must satisfy:*

(a) *Sparsity: with probability tending to one,* $\widehat{\boldsymbol{\beta}}(\text{ose})_2 = 0$.
(b) *Asymptotic normality:* $\sqrt{n}(\widehat{\boldsymbol{\beta}}(\text{ose})_1 - \boldsymbol{\beta}_{10}) \to N(0, I_1^{-1}(\boldsymbol{\beta}_{10}))$.

*In addition, consider* $p_{\lambda_n}(\cdot) = \lambda_n p(\cdot)$. *Suppose* $p'(\cdot)$ *is continuous on* $(0, \infty)$ *and there is some* $s > 0$ *such that* $p'(\theta) = O(\theta^{-s})$ *as* $\theta \to 0+$. *Then* (a) *and* (b) *hold, if* $n^{(1+s)/2}\lambda_n \to \infty$ *and* $\sqrt{n}\lambda_n \to 0$.

In Theorems 4 and 5, we have established the oracle properties of the one-step SCAD estimator. It is interesting to note that the choice of $\lambda_n$ is the same as that in Theorem 2 of Fan and Li [10]. It is also worth noting that our results require less regularity conditions than Theorem 2 of Fan and Li [10], for the penalty function does not need to be twice differentiable.

3.3. *Continuity of the one-step estimator.* For the nonconcave penalized likelihood estimates to be continuous, the minimum of the function $|\theta| + p'_\lambda(|\theta|)$ must be attained at 0 (Fan and Li [10]). Bridge penalty $(0 < q < 1)$ fails to satisfy the continuity condition, thus it is considered suboptimal (Fan and Li [10]). Our results require weaker conditions to ensure a continuous thresholding estimator. Note that $\widehat{\boldsymbol{\beta}}(\mathrm{ose})$ is obtained through an $\ell_1$ penalized criterion. Therefore, we only require $p'_\lambda(|\theta|)$ is continuous for $|\theta| > 0$ to ensure the continuity of $\widehat{\boldsymbol{\beta}}(\mathrm{ose})$. Theorem 4 and Theorem 5 indicate that Bridge penalty, $p_\lambda(|\theta|) = \lambda|\theta|^q$ for $0 < q < 1$, can be used in the one-step estimation scheme and their one-step estimates are continuous.

There is another interesting implication of the continuity of $\widehat{\boldsymbol{\beta}}(\mathrm{ose})$. Suppose two penalty functions have very similar derivatives, then we expect their one-step estimators are very close, too. To illustrate this point, we consider the limiting one-step estimator with the $L_q$ penalty when $q \to 0+$:

$$\boldsymbol{\beta}_q^{(1)} = \arg\min \tfrac{1}{2}(\boldsymbol{\beta} - \boldsymbol{\beta}^{(0)})^\mathsf{T}[-\nabla^2\ell(\beta^{(0)})](\boldsymbol{\beta} - \boldsymbol{\beta}^{(0)}) + n\sum_{j=1}^{p} \lambda q|\beta_j^{(0)}|^{q-1}|\beta_j|.$$

For each fixed $q$, we are interested in the whole profile of $\boldsymbol{\beta}_q^{(1)}$ as a function of $\lambda$. Thus we can consider $\lambda^* = \lambda q$ as the effective regularization parameter. On the other hand, suppose we consider the one-step estimator with the logarithm penalty, $p_\lambda(|\beta|) = \lambda \log|\beta|$,

$$\boldsymbol{\beta}_{\log}^{(1)} = \arg\min \tfrac{1}{2}(\boldsymbol{\beta} - \boldsymbol{\beta}^{(0)})^\mathsf{T}[-\nabla^2\ell(\beta^{(0)})](\boldsymbol{\beta} - \boldsymbol{\beta}^{(0)}) + n\sum_{j=1}^{p} \lambda|\beta_j^{(0)}|^{-1}|\beta_j|.$$

PROPOSITION 2. *If $q \to 0+$, then the profile of $\boldsymbol{\beta}_q^{(1)}$ converges to the profile of $\boldsymbol{\beta}_{\log}^{(1)}$ in the sense that* $\lim_{q\to 0+} \boldsymbol{\beta}_q^{(1)}(\lambda/q) = \boldsymbol{\beta}_{\log}^{(1)}(\lambda)$, $\forall \lambda > 0$.

We make a note that the convexity of the LLA is crucial for Proposition 2. We demonstrate the continuity property of the one-step estimator in linear regression models with an orthogonal design. As can be seen from Figure 2, in orthogonal design the $L_{0.01}$ penalty and the logarithm penalty are equivalent to some discontinuous thresholding rules, but their one-step estimators yield continuous thresholding rules. Moreover, the one-step $L_{0.01}$ estimator with $\lambda = 200$ is very similar to the one-step logarithm estimator with $\lambda = 2$, which shows us an illustrative example of Proposition 2. We also show the SCAD thresholding and its one-step version in Figure 2. They are both continuous and unbiased for large coefficients, but they are not identical.

**4. Implementation.** In this section we show that the LLA allows an efficient implementation of the one-step sparse estimator. The key is to notice that solving $\boldsymbol{\beta}^{(1)}$ is not much different from solving LASSO. Standard quadratic programming software can be used to solve LASSO. The shooting algorithm also works well (Fu [17] and Yuan and Lin [35]). Efron et al. [8] proposed an efficient path algorithm called LARS for computing the entire solution path of LASSO. See also the homotopy algorithm by Osborne, Presnell and Turlach [29]. The LARS algorithm is a major breakthrough in the development of the LASSO-type methods. Zou and Hastie [36] modified the LARS algorithm to compute the solution paths of the elastic net. Rosset and Zhu [31] generalized the LARS type algorithm to a class of optimization problems with a LASSO penalty. The LARS algorithm was used to simplify the computations in an empirical Bayes model for LASSO (Yuan and Lin [34]).

We adopt the LARS idea in our implementation. Write $\mu_i = \mathbf{x}_i^{\mathsf{T}}\boldsymbol{\beta}$ and $\ell_i = \ell_i(\mu_i, y_i)$. Observe that

$$-\nabla^2 \ell(\boldsymbol{\beta}^{(0)}) = \mathbf{X}^{\mathsf{T}}\mathbf{D}\mathbf{X},$$

where $\mathbf{D}$ is a $n \times n$ diagonal matrix with

$$\mathbf{D}_{ii} = -\frac{d^2\ell_i(\mu_i)}{d\mu_i^2}\Big|_{\hat{\mu}_i}, \qquad \hat{\mu}_i = \mathbf{x}_i^{\mathsf{T}}\boldsymbol{\beta}^{(0)}, i = 1, \ldots, n.$$

In linear regression models, $\mathbf{D}_{ii} = 2$. We separately discuss the algorithm for two types of concave penalties.

TYPE 1. $p_\lambda(t) = \lambda p(t)$ and $p'(t) > 0$ for all $t$. Bridge penalties and the logarithm penalty belong to this category which also covers many other penalties. We propose the following algorithm to compute the one-step estimator.

ALGORITHM 1.

Step 1. Create working data by

$$\mathbf{x}_{ij}^* = \sqrt{\mathbf{D}_{ii}}\mathbf{x}_{ij}/p'(|\beta_j^{(0)}|) \quad \text{and} \quad y_i^* = \sqrt{\mathbf{D}_{ii}}\hat{\mu}_i,$$
$$i = 1, 2, \ldots, n; j = 1, \ldots, p.$$

Step 2. Apply the LARS algorithm to solve

$$\hat{\boldsymbol{\beta}}^* = \arg\min_{\boldsymbol{\beta}} \left\{ \frac{1}{2}\sum_{i=1}^n (y_i^* - \mathbf{x}_i^{*\mathsf{T}}\boldsymbol{\beta})^2 + n\lambda \sum_{j=1}^p |\beta_j|. \right\}$$

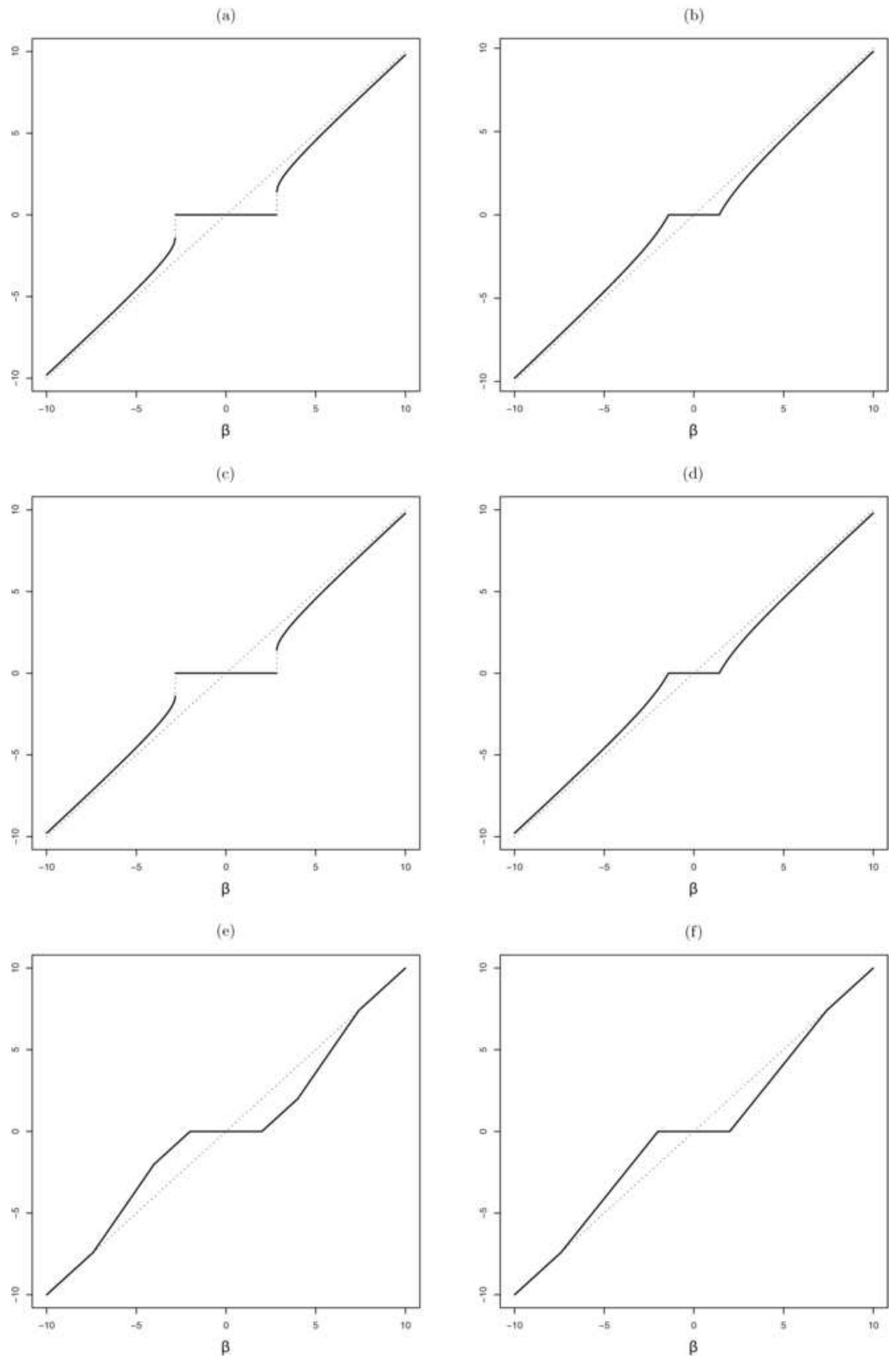

FIG. 2. *Compare thresholding rules in orthogonal design.* (a) *and* (b) *are for the logarithm penalty and its one-step LLA approximation,* $\lambda = 2$. (c) *and* (d) *are for Bridge* ($L_{0.01}$) *and its one-step LLA approximation,* $\lambda = 200$. (e) *and* (f) *are for SCAD and its one-step LLA approximation,* $\lambda = 2$.

Then, it is not hard to show that

$$\beta_j^{(1)} = \hat{\beta}_j^* / p'_\lambda(|\beta_j^{(0)}|), \qquad j = 1, 2, \ldots, p.$$

Thus, if $\hat{\beta}_j^* \neq 0$, then $x_j$ is selected in the final model.

Type 2. For some penalties, the derivative can be zero. In addition, the regularization parameter $\lambda$ cannot be separated from the penalty function. The SCAD penalty is a typical example. Let us assume that

$$U = \{j : p'_\lambda(|\beta_j^{(0)}|) = 0\} \quad \text{and} \quad V = \{j : p'_\lambda(|\beta_j^{(0)}|) > 0\}.$$

We write

$$\mathbf{X} = [\mathbf{X}_U, \mathbf{X}_V] \quad \text{and} \quad \boldsymbol{\beta}^{(1)} = ({\boldsymbol{\beta}^{(1)}}_U^{\mathsf{T}}, {\boldsymbol{\beta}^{(1)}}_V^{\mathsf{T}})^{\mathsf{T}}.$$

We propose the following algorithm to compute $\boldsymbol{\beta}^{(1)}$.

Algorithm 2.
Step 1a. Create working data by $y_i^* = \sqrt{\mathbf{D}_{ii}}\hat{\mu}_i$, $\mathbf{x}_i^* = \sqrt{\mathbf{D}_{ii}}\mathbf{x}_i$, $i = 1, \ldots, n$;
Step 1b. Let $\mathbf{x}_j^* = \mathbf{x}_j^* \frac{\lambda}{p'_\lambda(\beta_j^{(0)})}$ for $j \in V$.

Step 1c. Let $\mathbf{H}_U$ be the projection matrix in the space of $\{\mathbf{x}_j^*, j \in U\}$. Compute $y^{**} = y^* - \mathbf{H}_U y^*$ and $\mathbf{X}_V^{**} = \mathbf{X}_V^* - \mathbf{H}_U \mathbf{X}_V^*$.

Step 2. Apply the LARS algorithm to solve

$$\hat{\beta}_V^* = \arg\min_\beta \{\tfrac{1}{2}\|y^{**} - \mathbf{X}_V^{**}\boldsymbol{\beta}\|^2 + n\lambda\|\beta\|_1\}.$$

Step 3. Compute $\hat{\beta}_U^* = (\mathbf{X}_U^{*T}\mathbf{X}_U^*)^{-1}\mathbf{X}_U^{*T}(y^* - \mathbf{X}_V^*\hat{\beta}_V^*)$.

Then, it is not hard to show that

$$\beta_U^{(1)} = \hat{\beta}_U^* \quad \text{and} \quad \beta_j^{(1)} = \hat{\beta}_j^* \frac{\lambda}{p'_\lambda(|\beta_j^{(0)}|)} \qquad \text{for } j \in V.$$

Thus, if $\hat{\beta}_j^* \neq 0$, then $x_j$ is selected in the final model for $j \in V$.

In both algorithms the LARS step uses the same order of computations of a single OLS fit (Efron et al. [8]). Thus it is very efficient to compute the one-step estimator. It is also remarkable that if the penalty is of type 1, then the entire profile of the one-step estimator (as a function of $\lambda$) can be efficiently constructed. For the SCAD type penalty, we still need to solve the one-step estimator for each fixed $\lambda$, for the sets $U$ and $V$ could change as $\lambda$ varies.

**5. Numerical examples.** In this section we assess the finite sample performance of the one-step sparse estimates for linear regression models, logistic regression models and Poisson regression models in terms of model complexity (sparsity) and model error, defined by

$$ME\{\hat{\mu}(\cdot)\} = E\{\hat{\mu}(\mathbf{x}) - \mu(\mathbf{x})\}^2$$

for a selected model $\hat{\mu}(\cdot)$, where the expectation is taken over the new observation $\mathbf{x}$. We compare their performance with that of the SCAD with the original LQA algorithm (Fan and Li [10]) and the perturbed LQA algorithm (Hunter and Li [20]), and the best subset variable selection with the AIC, and BIC. For a fitted subset model $\mathcal{M}$, the AIC and BIC statistics are of the form

$$2\log(\text{likelihood}) - \lambda \cdot |\mathcal{M}|,$$

where $|\mathcal{M}|$ is the size of the model and $\lambda = 2$ and $\log(n)$, respectively. Note that the BIC is a consistent model selection criterion, while AIC is not. We further demonstrate the proposed methodology by analysis of a real data set.

In our simulation studies, we examine the performance of one-step sparse estimates with the SCAD penalty, logarithm penalty (defined in Section 3.3) and $L_{0.01}$ penalty. Note that we expect the logarithm penalty and $L_{0.01}$ penalty generate similar one-step sparse estimators. In Tables 1–3, one-step SCAD, one-step LOG and one-step $L_{0.01}$ stand for the one-step sparse estimate with the SCAD, logarithm and $L_{0.01}$ penalty, respectively; SCAD and P-SCAD represent the penalized least squares or likelihood estimators with the SCAD penalty using LQA and perturbed LQA algorithm, respectively; and AIC and BIC are the best subset variable selection with the AIC and BIC criteria, respectively. For the best subset variable selection, we exhaustively searched over all possible subsets. We used five-fold cross-validation to select the tuning parameters.

EXAMPLE 1 (*Linear model*). In this example, simulation data were generated from the linear regression model,

$$y = \mathbf{x}^{\mathsf{T}}\boldsymbol{\beta} + \epsilon,$$

where $\boldsymbol{\beta} = (3, 1.5, 0, 0, 2, 0, 0, 0, 0, 0, 0, 0)^{\mathsf{T}}$, $\epsilon \sim N(0, 1)$ and $\mathbf{x}$ is multivariate normal distribution with zero mean and covariance between the $i$th and $j$th elements being $\rho^{|i-j|}$ with $\rho = 0.5$. In our simulation, the sample size $n$ is set to be 50 and 100. For each case, we repeated the simulation 1,000 times.

For linear model, model error for $\hat{\mu} = \mathbf{x}^{\mathsf{T}}\widehat{\boldsymbol{\beta}}$ is $ME(\hat{\mu}) = (\widehat{\boldsymbol{\beta}} - \boldsymbol{\beta})^{\mathsf{T}} E(\mathbf{x}\mathbf{x}^{\mathsf{T}})(\widehat{\boldsymbol{\beta}} - \boldsymbol{\beta})$. Simulation results are summarized in Table 1, in which MRME stands for median of ratios of ME of a selected model to that of the ordinary least squares estimate under the full model. Both the columns of "C" and "IC" are measures of model complexity. Column "C" shows the average number of nonzero coefficients correctly estimated to be nonzero, and column "IC" presents the average number of zero coefficients incorrectly estimated to be nonzero. In the column labeled "Under-fit," we presented the proportion of excluding any nonzero coefficients in 1,000 replications. Likewise, we reported the probability of selecting the exact subset model and the probability of including all three significant variables and some noise variables in the columns "Correct-fit" and "Over-fit," respectively.

As can be seen from Table 1, all variable selection procedures dramatically reduce model error. One-step SCAD has the smallest model error among all competitors, followed by the SCAD and perturbed-SCAD. In terms of model error, penalized least squares methods with concave penalties outperform the best subset selection. In terms of sparsity, one-step SCAD also has the highest probability of correct fit. The SCAD penalty performs better than the other penalties. One-step LOG and one-step $L_{0.01}$ perform very similarly, which numerically confirms the assertion in Proposition 2. It is also interesting to note that a simulation study by Leng, Lin and Wahba [25] showed that in this example the LASSO did not consistently select the true model when optimizing the prediction error. In contrast, the noncon-cave penalty methods and their one-step estimates all work very well in this example because of their oracle properties.

TABLE 1
*Simulation results for linear regression models*

| | | No. of Zeros | | Proportion of | | |
|---|---|---|---|---|---|---|
| **Method** | **MRME** | **C** | **IC** | **Under-fit** | **Correct-fit** | **Over-fit** |
| | | | $n = 50$ | | | |
| One-step SCAD | 0.208 | 3.00 | 0.55 | 0.000 | 0.771 | 0.229 |
| One-step LOG | 0.263 | 3.00 | 0.89 | 0.000 | 0.559 | 0.441 |
| One-step $L_{0.01}$ | 0.262 | 3.00 | 0.90 | 0.000 | 0.555 | 0.445 |
| SCAD | 0.233 | 3.00 | 0.83 | 0.000 | 0.682 | 0.318 |
| P-SCAD | 0.235 | 3.00 | 0.64 | 0.000 | 0.701 | 0.299 |
| AIC | 0.660 | 3.00 | 1.84 | 0.000 | 0.195 | 0.805 |
| BIC | 0.401 | 3.00 | 0.63 | 0.000 | 0.576 | 0.424 |
| | | | $n = 100$ | | | |
| One-step SCAD | 0.234 | 3.00 | 0.55 | 0.000 | 0.784 | 0.216 |
| One-step LOG | 0.281 | 3.00 | 0.71 | 0.000 | 0.657 | 0.343 |
| One-step $L_{0.01}$ | 0.281 | 3.00 | 0.71 | 0.000 | 0.657 | 0.343 |
| SCAD | 0.252 | 3.00 | 0.75 | 0.000 | 0.732 | 0.268 |
| P-SCAD | 0.262 | 3.00 | 0.63 | 0.000 | 0.711 | 0.289 |
| AIC | 0.676 | 3.00 | 1.63 | 0.000 | 0.192 | 0.808 |
| BIC | 0.337 | 3.00 | 0.32 | 0.000 | 0.728 | 0.272 |

EXAMPLE 2 (*Logistic regression*). In this example, we simulated 1,000 data sets consisting of $n = 200$ observations from the model

$$Y|\mathbf{x} \sim \text{Bernoulli}\{p(\mathbf{x}^\mathsf{T}\boldsymbol{\beta})\},$$

where $p(u) = \exp(u)/(1 + \exp(u))$, and $\boldsymbol{\beta}$ is the same as that in Example 1. The covariate vector $\mathbf{x}$ is created as follows. We first generate $\mathbf{z}$ from a 12-dimensional multivariate normal distribution with zero mean and covariance between the $i$th and $j$th elements being $\rho^{|i-j|}$ with $\rho = 0.5$. Then we set $x_{2k-1} = z_{2k-1}$ and $x_{2k} = I(z_{2k} < 0)$ for $k = 1, \ldots, 6$, where $I(\cdot)$ is an indicator function. Thus, $\mathbf{x}$ has continuous as well as binary components.

Unlike the model error for linear regression models, there is no closed form of model error for the logistic regression model. In this example, the model error was estimated using Monte Carlo simulation. Simulation results are summarized in Table 2, in which MRME stands for median of ratios of ME of a selected model to that of the un-penalized maximum likelihood estimate under the full model, and other notation is the same as that in Table 1.

From Table 2, it can be seen that the best subset variable selection with the BIC criterion performs the best, however, the computational cost of the best subset variable selection is much more expensive than that of the nonconcave penalized likelihood approach. One-step sparse estimates require the least computational cost. It is interesting to see from Table 2 that the one-step SCAD performs as well as the fully iterative SCAD estimates by the LQA and perturbed LQA algorithms in terms of model error. The one-step estimates with logarithm and $L_{0.01}$ penalties perform very well. They have lower model error and rate of under-fit models than ones with the SCAD penalty.

TABLE 2
*Simulation results for logistic regression model*

| | | No. of Zeros | | Proportion of | | |
|---|---|---|---|---|---|---|
| **Method** | **MRME** | **C** | **IC** | **Under-fit** | **Correct-fit** | **Over-fit** |
| One-step SCAD | 0.238 | 2.95 | 0.82 | 0.051 | 0.565 | 0.384 |
| One-step LOG | 0.229 | 2.97 | 0.61 | 0.029 | 0.518 | 0.453 |
| One-step $L_{0.01}$ | 0.230 | 2.97 | 0.61 | 0.028 | 0.516 | 0.456 |
| SCAD | 0.238 | 2.92 | 0.51 | 0.076 | 0.706 | 0.218 |
| P-SCAD | 0.237 | 2.92 | 0.50 | 0.079 | 0.707 | 0.214 |
| AIC | 0.596 | 2.98 | 1.56 | 0.021 | 0.216 | 0.763 |
| BIC | 0.208 | 2.95 | 0.22 | 0.053 | 0.800 | 0.147 |

EXAMPLE 3 (*Poisson log-linear regression*). In this example, we considered a Poisson regression model

$$Y|\mathbf{x} \sim \text{Poisson}\{\lambda(\mathbf{x}^\mathsf{T}\boldsymbol{\beta})\},$$

where $\lambda(u) = \exp(u)$, $\boldsymbol{\beta} = (1.2, 0.6, 0, 0, 0.8, 0, 0, 0, 0, 0, 0, 0)^\mathsf{T}$ and $\mathbf{x}$ is the same as that of Example 1. We let the sample size be 60 and 120. For each case we simulated 1,000 data sets. Note that the model error is $ME(\widehat{\boldsymbol{\beta}}) = E\{\exp(\mathbf{x}^\mathsf{T}\widehat{\boldsymbol{\beta}}) - \exp(\mathbf{x}^\mathsf{T}\boldsymbol{\beta})\}^2$. Since $\mathbf{x}$ is normally distributed, we can derive a closed form for the model error using the moment generating function of normal distribution. Simulation results are summarized in Table 3, in which notation is the same as that in Table 2.

From Table 3, we can see that one-step SCAD sparse estimate outperforms the SCAD using both the original LQA algorithm and perturbed LQA algorithm in terms of model errors, model complexity and the rate of correct-fit. The best subset variable selection has the best rate of correct-fit for both $n = 60$ and 120. The correct-fit rate of one-step sparse estimates becomes much higher when the sample size increases from 60 to 120. This is not case for SCAD, P-SCAD and the best subset variable selection procedures.

EXAMPLE 4 (*Data analysis*). In this example, we demonstrate our one-step estimation methodology using the burns data, collected by the General

TABLE 3
*Simulation results for Poisson regression models*

| | | No. of Zeros | | Proportion of | | |
|---|---|---|---|---|---|---|
| **Method** | **MRME** | **C** | **IC** | **Under-fit** | **Correct-fit** | **Over-fit** |
| | | | $n = 60$ | | | |
| One-step SCAD | 0.284 | 2.99 | 1.35 | 0.011 | 0.386 | 0.603 |
| One-step LOG | 0.260 | 2.99 | 1.10 | 0.006 | 0.460 | 0.534 |
| One-step $L_{0.01}$ | 0.260 | 2.99 | 1.10 | 0.006 | 0.460 | 0.534 |
| SCAD | 0.292 | 3.00 | 2.75 | 0.003 | 0.095 | 0.902 |
| P-SCAD | 0.327 | 2.91 | 1.72 | 0.055 | 0.270 | 0.675 |
| AIC | 0.496 | 3.00 | 1.40 | 0.001 | 0.265 | 0.734 |
| BIC | 0.228 | 3.00 | 0.34 | 0.002 | 0.735 | 0.263 |
| | | | $n = 120$ | | | |
| One-step SCAD | 0.271 | 3.00 | 1.00 | 0.001 | 0.552 | 0.447 |
| One-step LOG | 0.266 | 3.00 | 0.76 | 0.000 | 0.603 | 0.397 |
| One-step $L_{0.01}$ | 0.266 | 3.00 | 0.77 | 0.000 | 0.601 | 0.399 |
| SCAD | 0.342 | 3.00 | 2.36 | 0.000 | 0.174 | 0.826 |
| P-SCAD | 0.356 | 2.95 | 1.60 | 0.037 | 0.322 | 0.641 |
| AIC | 0.594 | 3.00 | 1.45 | 0.000 | 0.235 | 0.765 |
| BIC | 0.277 | 3.00 | 0.25 | 0.000 | 0.790 | 0.210 |

Hospital Burn Center at the University of Southern California. The data set consists of 981 observations. Fan and Li [10] analyzed this data set as an illustration of the nonconcave penalized likelihood methods. As in Fan and Li [10], the binary response variable is taken to be the indicator whether the victims survived their burns or not. Four covariates, $x_1 = \text{age}$, $x_2 = \text{sex}$, $x_3 = \log(\text{burn area} + 1)$ and binary variable $x_4 = \text{oxygen}$ ($0 = \text{normal}$ $1 = \text{abnormal}$), are considered. To reduce modeling bias, quadratic terms of $x_1$ and $x_3$ and all interaction terms were included in the logistic regression model. We computed the one-step estimators with the SCAD and logarithm penalties. The regularization parameter was chosen by 5-fold cross-validation. The logarithm of selected $\lambda$ equals $-0.356$ and $-7.095$ for the one-step estimates with the SCAD and logarithm penalties, respectively.

With the selected regularization parameter, the fitted one-step SCAD sparse estimate yields the following model

$$\operatorname{logit}\{P(Y=1|\mathbf{x})\} = 4.82 - 8.74x_1 - 4.79x_3^2 + 6.67x_1x_3, \tag{5.1}$$

where $Y = 1$ stands for a victims survived from his/her burns. This model indicates that only $x_1$ and $x_3$ are significant. This is the same as the ones in the model selected by the SCAD with the LQA algorithm and reported in Fan and Li [10]. The one-step fit with logarithm penalty is

$$\begin{aligned}\operatorname{logit}\{P(Y=1|\mathbf{x})\} &= 4.55 - 6.45x_1 - 0.29x_4 \\ &\qquad - 0.56x_1^2 - 4.21x_3^2 + 5.21x_1x_3 - 0.15x_2x_3.\end{aligned} \tag{5.2}$$

It selects more variables than (5.1). This is consistent with Table 2, from which we can see that one-step fit with logarithm penalty has a higher rate of "over-fit" than the one-step SCAD estimator. The one-step $L_{0.01}$ fit is almost identical to (5.2).

## 6. Proofs.

6.1. *Proof of Theorem 1.* At the $k$-step, define a function with parameter $\boldsymbol{\beta}^{(k)}$ as follows

$$G(\boldsymbol{\beta}|\boldsymbol{\beta}^{(k)}) = \ell(\boldsymbol{\beta}) - n\sum_{j=1}^{p}[p_\lambda(|\beta_j^{(k)}|) + p'_\lambda(|\beta_j^{(k)}|)(|\beta_j| - |\beta_j^{(k)}|)].$$

Observe that $Q(\boldsymbol{\beta}^{(k)}) = G(\boldsymbol{\beta}^{(k)}|\boldsymbol{\beta}^{(k)})$, and

$$Q(\boldsymbol{\beta}) - G(\boldsymbol{\beta}|\boldsymbol{\beta}^{(k)}) = n\sum_{j=1}^{p}[p_\lambda(|\beta_j^{(k)}|) + p'_\lambda(|\beta_j^{(k)}|)(|\beta_j| - |\beta_j^{(k)}|) - p_\lambda(|\beta_j|)].$$

By the concavity of the penalty function $p_\lambda(\cdot)$, we have

$$p_\lambda(|\beta_j^{(k)}|) + p'_\lambda(|\beta_j^{(k)}|)(|\beta_j| - |\beta_j^{(k)}|) - p_\lambda(|\beta_j|) \geq 0.$$

If $\beta_j^{(k)} = 0$ we use the right derivative. Thus it follows that

$$Q(\boldsymbol{\beta}) \geq G(\boldsymbol{\beta}|\boldsymbol{\beta}^{(k)}).$$

We can take "$>$" in the above inequality if $p_\lambda(\cdot)$ is strictly concave. Moreover, it is easy to check that

$$\boldsymbol{\beta}^{(k+1)} = \arg\max_{\boldsymbol{\beta}} G(\boldsymbol{\beta}|\boldsymbol{\beta}^{(k)}).$$

Hence we have that

$$Q(\boldsymbol{\beta}^{(k+1)}) \geq G(\boldsymbol{\beta}^{(k+1)}|\boldsymbol{\beta}^{(k)}) \geq G(\boldsymbol{\beta}^{(k)}|\boldsymbol{\beta}^{(k)}) = Q(\boldsymbol{\beta}^{(k)}).$$

This completes the proof.

6.2. *Proof of Theorem 2.* Without loss of generality let us consider $t > t_0$. It suffices to show

$$\frac{\psi(t) - \psi^*(t)}{t - t_0} \geq 0. \tag{6.1}$$

Note that

$$\psi(t) - \psi^*(t) = \psi(t) - \psi(t_0) - p'_\lambda(t_0)(t - t_0).$$

Thus (6.1) is equivalent to

$$\frac{\psi(t) - \psi(t_0)}{t - t_0} \geq p'_\lambda(t_0). \tag{6.2}$$

Take a sequence of $\{t_k\}$ such that $t_0 < t_k < t$ and $t_k \to t_0$. By the convexity of $\phi(\cdot)$, we know

$$\frac{\psi(t) - \psi(t_0)}{t - t_0} \geq \frac{\psi(t_k) - \psi(t_0)}{t_k - t_0} \qquad \forall k. \tag{6.3}$$

Since $\phi(\cdot)$ is a majorization of $p_\lambda(\cdot)$ at $t_0$, we have

$$\frac{\psi(t_k) - \psi(t_0)}{t_k - t_0} \geq \frac{p_\lambda(t_k) - p_\lambda(t_0)}{t_k - t_0}. \tag{6.4}$$

Thus combining (6.3) and (6.4), we know

$$\frac{\psi(t) - \psi(t_0)}{t - t_0} \geq \frac{p_\lambda(t_k) - p_\lambda(t_0)}{t_k - t_0} \qquad \forall k.$$

Taking the limit in the above inequality we obtain (6.2). Similar arguments can be applied to the case of $t < t_0$.

6.3. *Proof of Theorem 3.* It suffices to show that

$$E\left[\frac{1}{\tau_j}\middle|\boldsymbol{\beta},\mathbf{y}\right]=np'_\lambda(|\beta_j|). \tag{6.5}$$

Then (2.13) is equivalent to (2.7), which in turn shows that LLA is identical to the EM algorithm.

By $p(\tau_j|\boldsymbol{\beta},\mathbf{y})\propto p(\beta_j|\tau_j)\pi(\tau_j)$, we have

$$E\left[\frac{1}{\tau_j}\middle|\boldsymbol{\beta},\mathbf{y}\right]=\frac{\int_0^\infty(1/\tau_j)p(\beta_j|\tau_j)\pi(\tau_j)\,d\tau_j}{\int_0^\infty p(\beta_j|\tau_j)\pi(\tau_j)\,d\tau_j},$$

and (2.11) and (2.12) yield

$$\frac{\int_0^\infty(1/\tau_j)p(\beta_j|\tau_j)\pi(\tau_j)\,d\tau_j}{\int_0^\infty p(\beta_j|\tau_j)\pi(\tau_j)\,d\tau_j}=-\frac{d\log(\exp(-np_\lambda(|\beta_j|)))}{d(|\beta_j|)}.$$

Hence (6.5) is proven.

By the nonnegativity of $H(t)$, it is easy to see that $\exp(-np_\lambda(|\beta|))$ is a strictly decreasing function of $|\beta|$, thus $p_\lambda(\cdot)$ is strictly increasing. To show $p_\lambda(\cdot)$ is unbounded, using dominant (or monotone) convergence theorem, we have $\exp(-np_\lambda(|\beta|))\to 0$ as $|\beta|\to\infty$. Hence $p_\lambda(\cdot)$ is unbounded.

6.4. *Proof of Theorem 4 and Theorem 5.* Theorem 4 can be proven by the same proof for Theorem 5, and therefore, we only prove Theorem 5.

Let us define

$$\begin{aligned}V_n(u)&=\frac{1}{2}\left(\frac{u}{\sqrt{n}}+\boldsymbol{\beta}_0-\boldsymbol{\beta}^{(0)}\right)^{\mathsf{T}}[-\nabla^2\ell(\beta^{(0)})]\left(\frac{u}{\sqrt{n}}+\boldsymbol{\beta}_0-\boldsymbol{\beta}^{(0)}\right)\\&\quad+n\sum_{j=1}^p p'_{\lambda_n}(|\beta_j^{(0)}|)\left|\beta_{0j}+\frac{u_j}{\sqrt{n}}\right|.\end{aligned}$$

$$\begin{aligned}V_n(u)-V_n(0)&=\frac{1}{2}\frac{u^{\mathsf{T}}}{\sqrt{n}}[-\nabla^2\ell(\beta^{(0)})]\frac{u}{\sqrt{n}}+(\boldsymbol{\beta}_0-\boldsymbol{\beta}^{(0)})^{\mathsf{T}}[-\nabla^2\ell(\beta^{(0)})]\frac{u}{\sqrt{n}}\\&\quad+n\sum_{j=1}^p p'_{\lambda_n}(|\beta_j^{(0)}|)\left(\left|\beta_{0j}+\frac{u_j}{\sqrt{n}}\right|-|\beta_{0j}|\right)\\&\equiv T_1+T_2+T_3.\end{aligned}$$

Let $\hat{u}(n)=\arg\min[V_n(u)-V_n(0)]$, then $\widehat{\boldsymbol{\beta}}(\mathrm{ose})=\boldsymbol{\beta}_0+\frac{\hat{u}(n)}{\sqrt{n}}$.

By Slutsky's theorem, it follows that

$$T_1=\frac{1}{2}\frac{u^{\mathsf{T}}}{\sqrt{n}}[-\nabla^2\ell(\beta^{(0)})]\frac{u}{\sqrt{n}}\xrightarrow{P}\frac{1}{2}u^{\mathsf{T}}I(\boldsymbol{\beta}_0)u. \tag{6.6}$$

$$T_2 = (\boldsymbol{\beta}_0 - \boldsymbol{\beta}^{(0)})^{\mathsf{T}}[-\nabla^2 \ell(\beta^{(0)})]\frac{u}{\sqrt{n}} \tag{6.7}$$
$$= \sqrt{n}(\boldsymbol{\beta}_0 - \boldsymbol{\beta}^{(0)})^{\mathsf{T}}\left[\frac{-\nabla^2 \ell(\beta^{(0)})}{n}\right]u \xrightarrow{D} -W^{\mathsf{T}} I(\boldsymbol{\beta}_0)u.$$

We can write $T_3$ as

$$T_3 = \sum_{j=1}^{p} \sqrt{n} p'_{\lambda_n}(|\beta_j^{(0)}|)\frac{|\beta_{0j} + u_j/\sqrt{n}| - |\beta_{0j}|}{1/\sqrt{n}} \equiv \sum_{j=1}^{p} T_{3j}.$$

Note that

$$\frac{|\beta_{0j} + u_j/\sqrt{n}| - |\beta_{0j}|}{1/\sqrt{n}} \to \mathrm{Sign}(\beta_{0j})u_j I(\beta_{0j} \neq 0) + |u_j| I(\beta_{0j} = 0).$$

We now examine the behavior of $\sqrt{n}p'_{\lambda_n}(|\beta_j^{(0)}|)$. First consider the case where $p'_{\lambda_n}(|\beta_j^{(0)}|) = \lambda_n p'(|\beta_j^{(0)}|)$. When $\beta_{0j} \neq 0$, since $|\beta_j^{(0)}| \to_P |\beta_{0j}|$, continuous mapping theorem says that $p'(|\beta_j^{(0)}|) \to_P p'(|\beta_{0j}|)$. Hence $\sqrt{n}\lambda_n \to 0$ yields $T_{3j} \to_P 0$. When $\beta_{0j} = 0$, $T_{3j} = 0$ if $u_j = 0$. For $u_j \neq 0$, we have

$$T_{3j} = |u_j|\sqrt{n}\lambda_n p'(|\beta_j^{(0)}|) = |u_j| n^{(1+s)/2}\lambda_n(|\sqrt{n}\beta_j^{(0)}|)^{-s}\frac{p'(|\beta_j^{(0)}|)}{|\beta_j^{(0)}|^{-s}}.$$

By $\sqrt{n}\beta_j^{(0)} \to_D N(0, I^{-1}(\boldsymbol{\beta}_0)_{jj})$, then from $n^{(1+s)/2}\lambda_n \to \infty$ we see $T_{3j} \to_P \infty$.

For the SCAD penalty, we have similar conclusions. $p'_\lambda(\theta) = 0$ if $\theta > a\lambda_n$ ($a = 3.7$). Thus, when $\beta_{0j} \neq 0$, $|\beta_j^{(0)}| \to_P |\beta_{0j}| > 0$, then $\lambda_n \to 0$ ensures $T_{3j} = \mathrm{Sign}(\beta_{0j})u_j\sqrt{n}p'_\lambda(|\beta_j^{(0)}|) \to_P 0$. When $\beta_{0j} = 0$, $T_{3j} = 0$ if $u_j = 0$. For $u_j \neq 0$, we have $|\beta_j^{(0)}| = O_p(\frac{1}{\sqrt{n}})$. Also note that $p'_\lambda(\theta) = \lambda_n$ for all $0 < \theta < \lambda_n$, which implies that if $\sqrt{n}\lambda_n \to \infty$, $T_{3j} = \sqrt{n}p'_\lambda(|\beta_j^{(0)}|)|u_j| = |u_j|\sqrt{n}\lambda_n$ with probability tending to one. Thus $T_{3j} \to_P \infty$.

Let us write $u = (u_{10}^{\mathsf{T}}, u_{20}^{\mathsf{T}})^{\mathsf{T}}$. Then we have

$$T_3 \to_P \begin{cases} 0, & \text{if } u_{20} = 0, \\ \infty, & \text{otherwise.} \end{cases} \tag{6.8}$$

Denote $W = (W_{10}^{\mathsf{T}}, W_{20}^{\mathsf{T}})^{\mathsf{T}}$. Combining (6.6), (6.7) and (6.8) we conclude that for each fixed $u$,

$$V_n(u) - V_n(0) \to_d V(u) \equiv \begin{cases} \frac{1}{2}u_{10}^{\mathsf{T}} I_1(\boldsymbol{\beta}_{10})u_{10} - W_{10}^{\mathsf{T}} u_{10}, & \text{if } u_{20} = 0, \\ \infty, & \text{otherwise.} \end{cases}$$

The unique minimum of $V(u)$ is $u_{10} = I_1^{-1}(\beta_{10})W_{10}$ and $u_{20} = 0$. $V_n(u) - V_n(0)$ is a convex function of $u$. By epiconvergence (Geyer [18] and Knight and Fu [21]), we conclude that

$$\hat{u}(n)_{10} \xrightarrow{d} I_1^{-1}(\beta_{10})W_{10}, \tag{6.9}$$

$$\hat{u}(n)_{20} \xrightarrow{d} 0. \tag{6.10}$$

By $W_{10} = N(0, I_1(\boldsymbol{\beta}_{10}))$, (6.9) is equivalent to

$$\sqrt{n}(\widehat{\boldsymbol{\beta}}(\text{ose})_1 - \boldsymbol{\beta}_{10}) \to N(0, I_1^{-1}(\boldsymbol{\beta}_{10})).$$

Note that (6.10) implies that $\sqrt{n}\widehat{\boldsymbol{\beta}}(\text{ose})_2 \to_P 0$. We now show that with probability tending to one, $\widehat{\boldsymbol{\beta}}(\text{ose})_2 = 0$. This is a stronger statement than (6.10). It suffices to prove that if $\beta_{0j} = 0$, $P(\hat{\beta}_j(\text{ose}) \neq 0) \to 0$. Assume $\hat{\beta}_j(\text{ose}) \neq 0$. By KKT conditions of (3.3), we must have

$$\frac{1}{\sqrt{n}}([-\nabla^2 \ell(\boldsymbol{\beta}^{(0)})](\widehat{\boldsymbol{\beta}}(\text{ose}) - \boldsymbol{\beta}^{(0)}))_j = \sqrt{n}\lambda_n p'_\lambda(|\beta_j^{(0)}|). \tag{6.11}$$

We have shown that when $\beta_{0j} = 0$, the right-hand side goes to $\infty$ in probability. However, the left-hand side can be written as

$$\left(\left[\frac{-\nabla^2 \ell(\boldsymbol{\beta}^{(0)})}{n}\right]\sqrt{n}(\widehat{\boldsymbol{\beta}}(\text{ose}) - \boldsymbol{\beta}_0)\right)_j - \left(\left[\frac{-\nabla^2 \ell(\boldsymbol{\beta}^{(0)})}{n}\right]\sqrt{n}(\widehat{\boldsymbol{\beta}}^{(0)} - \boldsymbol{\beta}_0)\right)_j.$$

By (6.9) and (6.10), we know the first term converges in law to some normal, and so does the second term. Thus

$$P(\hat{\beta}_j(\text{ose}) \neq 0) \leq P(\text{KKT condition}(6.11) \text{ holds}) \to 0.$$

**7. Discussion.** In this article, we have proposed a new algorithm based on the LLA for maximizing the nonconcave penalized likelihood. We further suggest using the one-step LLA estimator as the final estimates, because the one-step estimator naturally adopts a sparse representation and enjoys the oracle properties. In addition, the one-step sparse estimate can dramatically reduce the computational cost in the fully iterative methods. The simulation shows that one-step sparse estimates have very competitive performance with finite samples.

We have concentrated on the one-step sparse estimate for linear models and likelihood-based models, including generalized linear models. The proposed one-step sparse estimation method can be easily extended for variable selection in survival data analysis using penalized partial likelihood (Fan and Li [11] and Cai et al. [5]), variable selection for longitudinal data (Fan and Li [12]) and variable selection in semiparametric regression modeling (Li and Liang [26]).

**Acknowledgments.** The authors sincerely thank the co-editor, the associate editor and the referees for constructive comments that substantially improved an earlier version of this paper.

## REFERENCES

[1] Antoniadis, A. and Fan, J. (2001). Regularization of wavelets approximations. *J. Amer. Statist. Assoc.* **96** 939–967. MR1946364

[2] Bickel, P. J. (1975). One-step Huber estimates in the linear model. *J. Amer. Statist. Assoc.* **70** 428–434. MR0386168

[3] Black, A. and Zisserman, A. (1987). *Visual Reconstruction*. MIT Press, Cambridge, MA. MR0919733

[4] Breiman, L. (1996). Heuristics of instability and stabilization in model selection. *Ann. Statist.* **24** 2350–2383. MR1425957

[5] Cai, J., Fan, J., Li, R. and Zhou, H. (2005). Variable selection for multivariate failure time data. *Biometrika* **92** 303–316. MR2201361

[6] Cai, J., Fan, J., Zhou, H. and Zhou, Y. (2007). Marginal hazard models with varying-coefficients for multivariate failure time data. *Ann. Statist.* **35** 324–354.

[7] Cai, Z., Fan, J. and Li, R. (2000). Efficient estimation and inferences for varying-coefficient models. *J. Amer. Statist. Assoc.* **95** 888–902. MR1804446

[8] Efron, B., Hastie, T., Johnstone, I. and Tibshirani, R. (2004). Least angle regression. *Ann. Statist.* **32** 407–499. MR2060166

[9] Fan, J. and Chen, J. (1999). One-step local quasi-likelihood estimation. *J. Roy. Statist. Soc. Ser. B* **61** 927–943. MR1722248

[10] Fan, J. and Li, R. (2001). Variable selection via nonconcave penalized likelihood and its oracle properties. *J. Amer. Statist. Assoc.* **96** 1348–1360. MR1946581

[11] Fan, J. and Li, R. (2002). Variable selection for Cox's proportional hazards model and frailty model. *Ann. Statist.* **30** 74–99. MR1892656

[12] Fan, J. and Li, R. (2004). New estimation and model selection procedures for semiparametric modeling in longitudinal data analysis. *J. Amer. Statist. Assoc.* **99** 710–723. MR2090905

[13] Fan, J. and Li, R. (2006). Statistical challenges with high dimensionality: Feature selection in knowledge discovery. In *Proceedings of the Madrid International Congress of Mathematicians 2006* **3** 595–622. EMS, Zürich. MR2275698

[14] Fan, J., Lin, H. and Zhou, Y. (2006). Local partial likelihood estimation for life time data. *Ann. Statist.* **34** 290–325. MR2275243

[15] Fan, J. and Peng, H. (2004). On non-concave penalized likelihood with diverging number of parameters. *Ann. Statist.* **32** 928–961. MR2065194

[16] Frank, I. and Friedman, J. (1993). A statistical view of some chemometrics regression tools. *Technometrics* **35** 109–148.

[17] Fu, W. (1998). Penalized regression: The bridge versus the lasso. *J. Comput. Graph. Statist.* **7** 397–416. MR1646710

[18] Geyer, C. (1994). On the asymptotics of constrainted $M$-estimation. *Ann. Statist.* **22** 1993–2010. MR1329179

[19] Heiser, W. (1995). *Convergent Computation by Iterative Majorization: Theory and Applications in Multidimensional Data Analysis*. Clarendon Press, Oxford. MR1380319

[20] Hunter, D. and Li, R. (2005). Variable selection using mm algorithms. *Ann. Statist.* **33** 1617–1642. MR2166557

[21] Knight, K. and Fu, W. (2000). Asymptotics for lasso-type estimators. *Ann. Statist.* **28** 1356–1378. MR1805787

[22] Lange, K. (1995). A gradient algorithm locally equivalent to the EM algorithm. *J. Roy. Statist. Soc. Ser. B* **57** 425–437. MR1323348

[23] Lange, K., Hunter, D. and Yang, I. (2000). Optimization transfer using surrogate objective functions (with discussion). *J. Comput. Graph. Statist.* **9** 1–59. MR1819865

[24] Lehmann, E. and Casella, G. (1998). *Theory of Point Estimation*, 2nd ed. Springer, Berlin. MR1639875

[25] Leng, C., Lin, Y. and Wahba, G. (2006). A note on the lasso and related procedures in model selection. *Statist. Sinica* **16** 1273–1284.

[26] Li, R. and Liang, H. (2008). Variable selection in semiparametric regression modeling. *Ann. Statist.* **36** 261–286.

[27] Mike, W. (1984). Outlier models and prior distributions in Bayesian linear regression. *J. Roy. Statist. Soc. Ser. B* **46** 431–439. MR0790630

[28] Miller, A. (2002). *Subset Selection in Regression*, 2nd ed. Chapman and Hall, London. MR2001193

[29] Osborne, M., Presnell, B. and Turlach, B. (2000). A new approach to variable selection in least squares problems. *IMA J. Numer. Anal.* **20** 389–403. MR1773265

[30] Robinson, P. (1988). The stochastic difference between econometrics and statistics. *Econometrics* **56** 531–547. MR0946120

[31] Rosset, S. and Zhu, J. (2007). Piecewise linear regularized solution paths. *Ann. Statist.* **35** 1012–1030.

[32] Tibshirani, R. (1996). Regression shrinkage and selection via the lasso. *J. Roy. Statist. Soc. Ser. B* **58** 267–288. MR1379242

[33] Wu, Y. (2000). Optimization transfer using surrogate objective functions: Discussion. *J. Comput. Graph. Statist.* **9** 32–34. MR1819865

[34] Yuan, M. and Lin, Y. (2005). Efficient empirical Bayes variable selection and estimation in linear models. *J. Amer. Statist. Assoc.* **100** 1215–1225. MR2236436

[35] Yuan, M. and Lin, Y. (2006). Model selection and estimation in regression with grouped variables. *J. Roy. Statist. Soc. Ser. B* **68** 49–67. MR2212574

[36] Zou, H. and Hastie, T. (2005). Regularization and variable selection via the elastic net. *J. Roy. Statist. Soc. Ser. B* **67** 301–320. MR2137327

School of Statistics
University of Minnesota
Minneapolis, Minnesota 55455
USA
E-mail: hzou@stat.umn.edu

Department of Statistics
and The Methodology Center
Pennsylvania State University
University Park, Pennsylvania 16802–2111
USA
E-mail: rli@stat.psu.edu